\documentclass[a4paper,8pt]{article}
\usepackage[utf8x]{inputenc}
\usepackage{amsthm,amsmath,amssymb,oldgerm}
\usepackage[margin=3cm]{geometry}
\theoremstyle{plain}
\newtheorem{thm}{Theorem}
\newtheorem{lem}[thm]{Lemma}

\theoremstyle{definition}
 
\theoremstyle{definition}

\let\Im\relax
\DeclareMathOperator{\Im}{Im}
\DeclareMathOperator{\del}{\Delta_{hyp}}
\DeclareMathOperator{\deli}{\Delta_{hyp}^{{\it{i}}}}
\DeclareMathOperator{\delone}{\Delta_{hyp}^{1}}
\DeclareMathOperator{\deltwo}{\Delta_{hyp}^{2}}
\DeclareMathOperator{\hyp}{\mu_{hyp}}
\DeclareMathOperator{\hypvol}{\mu_{hyp}^{vol}}
\DeclareMathOperator{\hypone}{\mu_{hyp}^{1}}
\DeclareMathOperator{\hyptwo}{\mu_{hyp}^{2}}
\DeclareMathOperator{\hypi}{\mu_{hyp}^{{\it{i}}}}   
\DeclareMathOperator{\canvol}{\mu_{can}^{vol}}
\DeclareMathOperator{\canone}{{\it\mu_{\mathrm{can}}^{1}}}
\DeclareMathOperator{\cantwo}{{\it{\mu_{\mathrm{can}}^{2}}}}
\DeclareMathOperator{\cani}{\mu_{can}^{{\it{i}}}}
\DeclareMathOperator{\hatcani}{\overline{\mu}_{can}^{{\it{i}}}}
\DeclareMathOperator{\khypi}{{\it{K_{\mathrm{hyp}}^{{\it{i}}}}}}
\DeclareMathOperator{\khypone}{{\it{K_{\mathrm{hyp}}^{1}}}}
\DeclareMathOperator{\khyptwo}{{\it{K_{\mathrm{hyp}}^{2}}}}

\DeclareMathOperator{\vx}{\mathrm{vol_{\mathrm{hyp}}}}

\DeclareMathAlphabet{\mathpzc}{OT1}{pzc}{m}{it}
\makeatletter

\newcommand{\Rmnum}[1]{\expandafter\@slowromancap\romannumeral #1@}
\title{A relation of cusp forms and Maass forms on product of hyperbolic Riemann orbisurfaces of finite volume}
\author{Anilatmaja Aryasomayajula}
\date{}
\begin{document}
\maketitle
\begin{abstract}
\noindent In \cite{jk}, J.~Jorgenson and J.~Kramer proved a certain key identity which relates the two natural metrics, 
namely the hyperbolic metric and the canonical metric defined on a compact hyperbolic Riemann surface. 
In this article, we extend this identity to product of noncompact hyperbolic Riemann orbisurfaces of 
finite volume, which can be realized as a quotient space of the action of a Fuchsian subgroup of first kind on the hyperbolic upper half plane. 

\vspace{0.2cm}\noindent
\end{abstract}
\section{Introduction}
For $i=1,2$, let $X_{i}$ be a noncompact hyperbolic Riemann orbisurface of finite volume $\vx(X_{i})$ 
with genus $g_{i}\geq 1$, and can be realized as the quotient space $\Gamma_{i}\backslash\mathbb{H}$, 
where $\Gamma_{i}\subset \mathrm{PSL}_{2}(\mathbb{R})$ is a Fuchsian subgroup of the first kind acting 
on the hyperbolic upper half-plane $\mathbb{H}$, via fractional linear transformations. Let $\mathcal{
E}_{i}$ and $\mathcal{P}_{i}$ denote the finite sets of elliptic fixed points and cusps of $\Gamma_{i}$, 
respectively. Put $\overline{X}_{i}=X_{i}\cup \mathcal{P}_{i}$. Then, $\overline{X}_{i}$ admits the structure of 
a Riemann surface. Now consider the complex surface $X=X_{1}\times X_{2}$, which admits the structure of a 
noncompact K\"ahler-orbifold of dimension two. Put $\overline{X}=\overline{X}_{1}\times \overline{X}_{2}$. 

\vspace{0.2cm}
For $i=1,2$, let $\hypi$ denote the (1,1)-form associated to hyperbolic metric, which is 
the natural metric on $X_{i}$, and of constant negative curvature minus one. Now $\hypone +\hyptwo$ 
is the natural metric on $X$, and let $\hypvol$ denote the volume form associated to 
$\hypone +\hyptwo$. 

\vspace{0.2cm}
For $i=1,2$, the Riemann surface $\overline{X}_{i}$ is embedded in its Jacobian variety $\mathrm{Jac}(
\overline{X}_{i})$ via the Abel-Jacobi map. Then, the pull back of the flat Euclidean metric by the Abel-Jacobi map is called the 
canonical metric, and the (1,1)-form associated to it is denoted by $\hatcani$. We denote its restriction 
to $X_{i}$ by $\cani$. Now, $\canone+\cantwo$ defines a metric on $X$, which corresponds to the flat 
Euclidean metric, and let $\canvol$ denote the volume form associated to $\canone+\cantwo$. 

\vspace{0.2cm}
For $i=1,2$, let $\deli$ denote the hyperbolic Laplacian acting on smooth functions on $X_{i}$. The 
hyperbolic heat kernel $\khypi(t;z_{i},w_{i})$ on $\mathbb{R}_{>0}\times X_{i}\times X_{i}$ is the unique 
solution of the heat equation 
\begin{align*}
\bigg(\deli + \frac{\partial}{\partial t}\bigg)\khypi(t;z_{i},w_{i}) =0,
\end{align*}
with the normalization condition
\begin{align*}  
\lim_{t\rightarrow 0}\int_{X_{i}}\khypi(t;z_{i},w_{i})f(z_{i})\hypi(z_{i}) = f(w_{i}),
\end{align*}
for any fixed $w\in X_{i}$ and any smooth function $f$ on $X_{i}$. When $z_{i}=w_{i}$, for brevity of 
notation, we denote the hyperbolic heat kernel by $\khypi(t;z_{i})$.

\vspace{0.3cm}\textbf{Main result}
With notation as above, for $z=(z_{1},z_{2})\in (X_{1}\backslash \mathcal{E}_{1})\times 
(X_{2}\backslash \mathcal{E}_{2})$, we have the relation of differential forms
\begin{align*}
&g_{1}g_{2}\canvol(z)=\bigg(\frac{1}{4\pi}+ \frac{1}{\vx(X_{1})}\bigg)\cdot \bigg(\frac{1}{4\pi}+ 
\frac{1}{\vx(X_{2})}\bigg)\hypvol(z)+\frac{1}{2}\bigg(\frac{1}{4\pi}+ \frac{1}{\vx(X_{1})}\bigg)\times
\\&\bigg(\int_{0}^{\infty}\deltwo \khyptwo(t;z_{2}) dt\bigg)\hypvol(z)+
\frac{1}{2}\bigg(\frac{1}{4\pi}+ \frac{1}{\vx(X_{2})}\bigg)\cdot\bigg(\int_{0}^{\infty}\delone \khypone(t;z_{1}) dt\bigg)\hypvol(z)+\notag\\&
\frac{1}{4}\bigg(\int_{0}^{\infty}\delone \khypone(t;z_{1}) dt\bigg)\cdot\bigg(\int_{0}^{\infty}\deltwo \khyptwo(t;z_{2}) dt\bigg)\hypvol(z).
\end{align*}
The above relation, which we call the key-identity, relates the two natural metrics defined on a 
K\"ahler-orbifold of dimension two. The key-identity is proved for compact hyperbolic Riemann surfaces, 
by J.~Jorgenson and J.~Kramer in \cite{jk}. The same authors extended the key-identity to noncompact 
hyperbolic Riemann surfaces of finite volume in \cite{bonn}. Following diferent methods, the 
key-identity is extended to noncompact hyperbolic Riemann orbisurfaces in \cite{anilpaper}.
\paragraph{Arithmetic significance}
The key-identity is the most significant technical result of \cite{jk}, which transforms a 
problem in Arakelov theory into that of hyperbolic geometry. The key-identity has enabled 
J.~Jorgenson and J.~Kramer to derive optimal bounds for the canonical Green's function defined on a compact hyperbolic Riemann 
surface $X$ in terms of invariants coming from the hyperbolic geometry of $X$.

\vspace{0.2cm}\noindent
Using the key-identity one can relate the holomorphic world of cusp forms with the $C^{\infty}$ 
world of M\"ass forms, via the spectral expansion of the hyperbolic heat kernel in terms 
of M\"ass forms. In fact, J.~Jorgenson and J.~Kramer have derived a Rankin-Selberg $L$-function 
relation relating the Fourier coefficients of cusp forms with those of M\"ass forms in \cite{bonn}.

\vspace{0.2cm}\noindent
Our main result is the first instance of an extension of the key-identity to higher dimensions, 
relating the cusp forms of the group $\Gamma_{1}\times\Gamma_{2}$ with the M\"ass forms defined on 
$X_{1}$ and $X_{2}$ via the spectral expansions of the hyperbolic heat kernels $\khypone(t;z_{1})$ and 
$\khyptwo(t;z_{2})$, respectively. 
{\small{\paragraph{Acknowledgements}
This article was realized during the author's graduate studies, which was completed under the 
supervision of  J.~Kramer at Humboldt Universit\"at zu Berlin. The author would like to express his 
gratitude to J.~Jorgenson and J.~Kramer for introducing him to the exciting area of heat kernels and 
automorphic forms, and for generously discussing many new scientific ideas, which resulted in the 
completion of this project. 

The author would also like to extend his gratitude to N.~Walji, who through many interesting scientific 
discussions has conveyed the importance of the key-identity. }} 
\section{Background material}\label{section1}
For $i=1, 2$, let $\Gamma_{i} \subset \mathrm{PSL}_{2}(\mathbb{R})$ be a Fuchsian subgroup of the first kind acting by 
fractional linear transformations on the upper half-plane $\mathbb{H}$. Let $X_{i}$ be the quotient space $\Gamma_{i}\backslash \mathbb{H}$, and let $g_{i}$ 
denote the genus of $X_{i}$. The quotient space $X_{i}$ admits the structure of a Riemann orbisurface.  

\vspace{0.2cm}
Let $\mathcal{E}_{i}$ and $\mathcal{P}_{i}$ be the finite sets of elliptic fixed points and cusps of 
$X_{i}$, respectively; put $\mathcal{S}_{i} = \mathcal{E}_{i}\cup\mathcal{P}_{i}$. For 
$\mathfrak{e}_{i}\in \mathcal{E}_{i}$, let $m_{\mathfrak{e}_{i}}$ denote the order of $\mathfrak{e}_{i}$; for 
$p_{i}\in \mathcal{P}_{i}$, put $m_{p_{i}}=\infty$; for $z_{i}\in X_{i}\backslash \mathcal{E}_{i}$, put $m_{z_{i}}=1.$ Let $\overline{X}_{i}$ 
denote $\overline{X}_{i}=X_{i}\cup\mathcal{P}_{i}.$ 

\vspace{0.2cm}
Locally, away from the elliptic fixed points and the cusps, we identify $\overline{X}_{i}$ with its 
universal cover $\mathbb{H}$, and hence, denote the points on $\overline{X}_{i}\backslash 
\mathcal{S}_{i}$ by the same letter as the points on $\mathbb{H}$.

\vspace{0.2cm}
The quotient space $\overline{X}_{i}$ admits the structure of a compact Riemann surface. We refer the 
reader to Section 1.8 in $\cite{miyake}$, for the details regarding the structure of $\overline{X}_{i}$ as a compact Riemann surface. 

\vspace{0.2cm}
Let $X$ denote the product of the Riemann orbisurfaces $X_{1}\times X_{2}$. Then, $X$ admits the 
structure of a complex K\"ahler-orbifold of dimension two. The boundary of $X$ is given by 
$\partial X=\big(X_{1}\times\mathcal{P}_{2}\big)\cup\big(X_{2}\times\mathcal{P}_{1}
\big)$, and the compactification of $X$ is given by $\overline{X}=\overline{X}_{1}\times\overline{X}_{2}$. 

\paragraph{Hyperbolic metric} 
For $i=1,2$, we denote the (1,1)-form corresponding to the hyperbolic metric of $X_{i}$, which is compatible with the complex structure on $X_{i}$ 
and has constant negative curvature equal to minus one, by $\hypi(z_{i})$. Locally, for 
$z_{i}\in X_{i}\backslash \mathcal{E}_{i}$, it is given by
\begin{equation*}
 \hypi(z_{i})= \frac{i}{2}\cdot\frac{dz_{i}\wedge d\overline{z}_{i}}{{\Im(z_{i})}^{2}}.
\end{equation*} 
From the above formula, it follows that the hyperbolic metric $\hypi(z_{i})$ is singular at the elliptic fixed points and at the cusps. 

\vspace{0.2cm}
Let $\vx(X_{i})$ be the volume of $X_{i}$ with respect to the hyperbolic metric $\hypi(z_{i})$. It is given by the formula 
\begin{equation*}
\vx(X_{i}) = 2\pi\bigg(2g_{i} -2 + |\mathcal{P}_{i}|+\sum_{\mathfrak{e}_{i} \in \mathcal{E}_{i}}\bigg(1-
\frac{1}{m_{\mathfrak{e}_{i}}}\bigg)\bigg). 
\end{equation*}
We denote the (1,1)-form corresponding to the hyperbolic metric of $X$, which is compatible with the complex structure on $X$, by $\hyp(z)$ 
and the corresponding volume form by $\hypvol(z)$. Locally, for 
$z=(z_{1},z_{2})\in (X_{1}\backslash \mathcal{E}_{1})\times(X_{2}\backslash \mathcal{E}_{2})$, it is given by
\begin{align*}
\hyp(z)= \hypone(z_{1})+\hyptwo(z_{2})=\frac{i}{2}\cdot\frac{dz_{1}\wedge d\overline{z}_{1}}{{\Im(z_{1})}^{2}}+\frac{i}{2}\cdot\frac{dz_{2}\wedge 
d\overline{z}_{2}}{{\Im(z_{2})}^{2}}, 
\end{align*}
and the corresponding volume form is given by 
\begin{align*}
\hypvol(z)=\hypone(z_{i})\wedge\hyptwo(z_{2})=-\frac{dz_{1}\wedge d\overline{z}_{1}\wedge dz_{2}\wedge d\overline{z}_{2}}{2\Im(z_{1})^{2}\Im(z_{2})^{2}}. 
\end{align*}
\paragraph{Canonical metric}
For $i=1,2$, let $S_{2}(\Gamma_{i})$ denote the $\mathbb{C}$-vector space of cusp forms of 
weight 2 with respect to $\Gamma_{i}$ equipped with the Petersson inner-product. Let 
$\lbrace f_{1}^{i},\ldots,f_{g_{i}}^{i}\rbrace $ denote an orthonormal basis of $S_{2}(\Gamma_{i})$ 
with respect to the Petersson inner-product. Then, the (1,1)-form $\cani(z_{i})$ corresponding to the 
canonical metric of $X_{i}$ is given by 
\begin{equation*}
\cani(z_{i})=\frac{i}{2g_{i}} \sum_{j=1}^{g_{i}}\left|f_{j}^{i}(z_{i})\right|^{2}dz_{i}\wedge d\overline{z}_{i}.
\end{equation*}
The canonical metric $\cani(z)$ remains smooth at the elliptic fixed points and at the cusps, and measures 
the volume of $X$ to be one. 

\vspace{0.2cm}
Let $\Omega_{\overline{X}_{i}}^{1}$ denote the cotangent bundle of holomorphic differential forms of degree 
one on $\overline{X}_{i}$. Recall that for each $f^{i}\in S_{2}(\Gamma_{i})$, $f^{i}(z_{i})dz_{i}$ 
defines a holomorphic differential form of degree one on $\overline{X}_{i}$, and every holomorphic 
differential form of degree one on $\overline{X}_{i}$ comes from a weight 2 cusp form. So 
$\lbrace f_{1}^{i},\ldots,f_{g_{i}}^{i}\rbrace $ the orthonormal basis of $S_{2}(\Gamma_{i})$ with 
respect to the Petersson inner-product gives us an orthonormal basis $\lbrace f_{1}^{i}dz_{i},\ldots,f_{g_{i}}^{i}dz_{i}\rbrace $ of $H^{0}(\overline{X}_{i},\Omega_{\overline{X}_{i}}^{1})$ 
endowed with the $L^{2}$-inner product given by 
\begin{align*}
\langle\alpha^{i},\beta^{i} \rangle=\frac{i}{2}\int_{\overline{X}}\alpha^{i}(z_{i})\overline{\beta^{i}(z_{i})}, 
\end{align*}
where $\alpha^{i},\beta^{i}\in\Omega_{\overline{X}_{i}}^{1}$.  

\vspace{0.2cm}
Let $\Omega_{\overline{X}}^{2}$ denote the space of holomorphic differential forms of degree 2, and let $\lbrace \omega_{1}, \ldots, \omega_{n} 
\rbrace $ denote an orthonormal basis of $H^{0}(\overline{X},\Omega_{\overline{X}}^{2})$ endowed with the $L^{2}$-inner product given by 
\begin{align*}
\langle\alpha,\beta \rangle=-\frac{1}{4}\int_{\overline{X}}\alpha(z)\overline{\beta(z)}, 
\end{align*}
where $n$ denotes the dimension of $H^{0}(\overline{X},\Omega_{\overline{X}}^{2})$ as a vector space over $\mathbb{C}$, and 
$\alpha,\beta\in\Omega_{\overline{X}}^{2}$. Then, the canonical volume form on $X$ is defined as
\begin{align*}
\canvol(z)=-\frac{1}{4n}\sum_{j=1}^{n}\omega_{j}(z)\wedge\overline{\omega_{j}(z)}.
\end{align*}
The canonical volume form $\canvol(z)$ measures the volume of $X$ to be one.    
\paragraph{Hyperbolic Laplacian}
For $i=1,2$, the hyperbolic Laplacian acting on smooth functions defined on $X_{i}$ is given by 
\begin{align*}
\deli =-y_{i}^{2}\bigg(\frac{\partial^{2}}{\partial x_{i}^{2}}+\frac{\partial^{2}}{\partial y_{i}^{2}}\bigg), 
\end{align*}
and the hyperbolic Laplacian acting on smooth functions defined on $X$ is given by 
\begin{align*}
\del=\delone+\deltwo= -y_{1}^{2}\bigg(\frac{\partial^{2}}{\partial x_{1}^{2}}+\frac{\partial^{2}}{\partial y_{1}^{2}}\bigg)-
y_{2}^{2}\bigg(\frac{\partial^{2}}{\partial x_{2}^{2}}+\frac{\partial^{2}}{\partial y_{2}^{2}}\bigg).
\end{align*}
\paragraph{Hyperbolic heat kernels}
For $t \in \mathbb{R}_{> 0}$  and $z, w \in \mathbb{H}$, let $K_{\mathbb{H}}(t;z,w)$ denote the hyperbolic heat kernel on $\mathbb{R}_{> 0}
\times\mathbb{H}\times\mathbb{H}$. 

\vspace{0.2cm}
For $i=1,2$, $t \in  \mathbb{R}_{> 0}$ and $z_{i}, w_{i} \in X_{i}$, the hyperbolic heat kernel $\khypi(t;z_{i},w_{i})$ on 
$\mathbb{R}_{> 0}\times X_{i}\times X_{i}$ is defined as 
\begin{equation}\label{defnkhyp}
\khypi(t;z_{i},w_{i})=\sum_{\gamma_{i}\in\Gamma_{i}}K_{\mathbb{H}}(t;z_{i},\gamma_{i} w_{i}).
\end{equation}
For $z_{i},w_{i}\in X_{i}$, the hyperbolic heat kernel $\khypi(t;z_{i},w_{i})$ satisfies the differential equation
\begin{align}\label{diffeqnkhyp} 
\bigg(\Delta_{\text{hyp},z}^{i} + \frac{\partial}{\partial t}\bigg)\khypi(t;z,w) &=0,
\end{align}
Furthermore, for a fixed $w_{i}\in X_{i}$, and any smooth function $f^{i}$ on $X_{i}$, the hyperbolic heat kernel $\khypi(t;z_{i},w_{i})$ 
satisfies the equation
\begin{align}\label{normcondkhyp}  
\lim_{t\rightarrow 0}\int_{X_{i}}\khypi(t;z_{i},w_{i})f^{i}(z)\hypi(z_{i}) &= f^{i}(w_{i}).
\end{align}
To simplify notation, we write $\khypi(t;z_{i})$ instead of $\khypi(t;z_{i},z_{i})$, when $z_{i}=w_{i}$.
\paragraph{Key-identity on $X_{i}$}
For $i=1,2$, $z_{i}\in X_{i} \backslash \mathcal{E}_{i}$, we have the relation of differential forms
\begin{align}\label{keyidentityrs}
&g_{i}\cani(z_{i}) =\bigg(\frac{1}{4\pi}+ \frac{1}{\vx(X_{i})}\bigg)\hypi(z_{i})+ \frac{1}{2}\bigg(\int_{0}^{\infty}\deli \khypi(t;z_{i}) dt\bigg)
\hypi(z_{i}).
\end{align}
This relation has been established as Theorem 3.4 in \cite{jk}, when $X_{i}$ is 
compact. The proof given in \cite{jk} applies to our case where $X$ does admit elliptic fixed points and 
cusps, as long as $z_{i}\in X_{i}\backslash \mathcal{E}_{i}$. 

\vspace{0.2cm}
In \cite{anilpaper}, the above identity is extended to elliptic fixed points and cusps at the level of 
currents. 
\section{Key-identity on $X$}\label{section2}
\begin{lem}\label{lem1}
The dimension of $H^{0}(\overline{X},\Omega_{\overline{X}}^{2})$ as a vector space over $\mathbb{C}$ is $g_{1}g_{2}$, and 
for $z=(z_{1},z_{2})\in X$, the canonical volume form is given by
\begin{align}\label{lem1eqn}
\canvol(z)=\canone(z_{1})\wedge\cantwo(z_{2})=-\frac{1}{4g_{1}g_{2}}\sum_{j=1}^{g_{1}}\sum_{k=1}^{g_{2}}|f^{1}_{j}(z_{1})|^{2}\cdot|f^{2}_{k}(z_{2})|^{2} 
dz_{1}\wedge d\overline{z}_{1}\wedge dz_{2}\wedge d\overline{z}_{2}.
\end{align}
\begin{proof}
From K\"unneth theorem of algebraic geometry, we have
\begin{align*}
H^{0}\big(\overline{X},\Omega_{\overline{X}}^{2}\big)=H^{0}\big(\overline{X}_{1},\Omega_{\overline{X}_{1}}^{1}\big)\otimes H^{0}\big(\overline{X}_{2},
\Omega_{\overline{X}_{2}}^{1}\big).
\end{align*}
So from the isomorphism $S_{2}(\Gamma_{i})\cong H^{0}(\overline{X}_{i},\Omega_{\overline{X}_{i}}^{1})$, it follows that the 
set 
\begin{align*}
\big\lbrace f^{1}_{j}f^{2}_{k} \big\rbrace_{\substack{1\leq j\leq g_{1}\\1\leq k\leq g_{2}}}
\end{align*}
forms an orthonormal basis of $H^{0}(\overline{X},\Omega_{\overline{X}}^{2})$, which implies that 
\begin{align*}
\canvol(z)=-\frac{1}{4g_{1}g_{2}}\sum_{j=1}^{g_{1}}\sum_{k=1}^{g_{2}}|f^{1}_{j}(z_{1})|^{2}\cdot |f^{2}_{k}(z_{2})|^{2} 
dz_{1}\wedge d\overline{z}_{1}\wedge dz_{2}\wedge d\overline{z}_{2}.
\end{align*}
Furthermore, a direct calculation, shows that $\canone(z_{1})\wedge\cantwo (z_{2})=\canvol(z)$, which completes the proof of the lemma. 
\end{proof}
\end{lem}
\begin{thm}\label{thm2}
For $z=(z_{1},z_{2})\in (X_{1}\backslash \mathcal{E}_{1})\times (X_{2}\backslash \mathcal{E}_{2})$, we have the relation of differential forms
\begin{align*}
&g_{1}g_{2}\canvol(z)=\bigg(\frac{1}{4\pi}+ \frac{1}{\vx(X_{1})}\bigg)\cdot \bigg(\frac{1}{4\pi}+ 
\frac{1}{\vx(X_{2})}\bigg)\hypvol(z)+\frac{1}{2}\bigg(\frac{1}{4\pi}+ \frac{1}{\vx(X_{1})}\bigg)\times
\\&\bigg(\int_{0}^{\infty}\deltwo \khyptwo(t;z_{2}) dt\bigg)\hypvol(z)+
\frac{1}{2}\bigg(\frac{1}{4\pi}+ \frac{1}{\vx(X_{2})}\bigg)\cdot\bigg(\int_{0}^{\infty}\delone \khypone(t;z_{1}) dt\bigg)\hypvol(z)+\notag\\&
\frac{1}{4}\bigg(\int_{0}^{\infty}\delone \khypone(t;z_{1}) dt\bigg)\cdot\bigg(\int_{0}^{\infty}\deltwo \khyptwo(t;z_{2}) dt\bigg)\hypvol(z).
\end{align*}
\begin{proof}
The proof of the theorem follows from combining equations \eqref{keyidentityrs} and \eqref{lem1eqn}.
\end{proof}
\end{thm}

\vspace{0.3cm}
{\small{
Department of Mathematics, \\University of Hyderabad, \\Prof. C.~R.~Rao Road, Gachibowli,\\
Hyderabad, 500046, India\\email: anilatmaja@gmail.com}}

\begin{thebibliography}{99}
\addcontentsline{toc}{chapter}{Bibliography}
\bibitem[1]{anilpaper} Anilatmaja Aryasomayajula, Extension of a key identity, arXiv:1310.4336.
\bibitem[2]{jk} J. Jorgenson and J. Kramer, Bounds on canonical Green's functions, Compositio Math. 
142 (2006), 679--700.
\bibitem[3]{bonn} J. Jorgenson and J. Kramer, A relation involving Rankin-Selberg L-functions 
of cusp forms and Maass forms, In: B. Krötz, O. Offen, E. Sayag (eds.), Representation Theory, 
Complex Analysis, and Integral Geometry, 9--40, Birkhäuser-Verlag, 2012.
\bibitem[4]{miyake} T. Miyake, Modular Forms, Springer-Verlag, Berlin, 2006.
\end{thebibliography}
\end{document}